\documentclass[aop,authoryear, preprint]{imsart} 

\usepackage{amsfonts,amsthm,amsmath,amssymb,amstext,graphicx}
\usepackage{mathrsfs}
\RequirePackage[OT1]{fontenc}
\RequirePackage{natbib}
\RequirePackage[dvips]{hyperref}
\RequirePackage{hypernat}

\startlocaldefs
\newcommand{\msf}[1]{\mathsf{#1}}
\newcommand{\Pb}{{\msf{P}}}
\newcommand{\Eb}{{\msf{E}}}
\newcommand{\Qb}{{\msf{Q}}}
\newcommand{\Fc}{{ \mathscr{F}}}
\newcommand{\Mc}{{ \mathscr{M}}}

\newcommand{\set}[1]{\left\{#1\right\}}

\newcommand{\brc}[1]{\left(#1\right)}

\renewcommand{\le}{\leqslant} 
\renewcommand{\ge}{\geqslant}

\newtheorem*{theorem}{Theorem}

\theoremstyle{definition} 

\endlocaldefs

\begin{document}

\begin{frontmatter}
\noindent {\sl{June 4, 2010}}
\title{On the Expectation of the First Exit Time of a Nonnegative Markov Process Started at a Quasistationary Distribution}
\runtitle{Expectation of the First Exit Time}


\begin{aug}

\author{\fnms{Moshe} \snm{Pollak}\ead[label=e1]{msmp@mscc.huji.ac.il}\thanksref{t1}}
\affiliation{Hebrew University of Jerusalem}
\address{M. Pollak \\
Hebrew University of Jerusalem\\
Department of Statistics\\
Mount Scopus, Jerusalem, 91905, Israel\\
\printead{e1}}
 \and
\author{\fnms{Alexander G.} \snm{Tartakovsky}\corref{}\ead[label=e2]{tartakov@usc.edu}\thanksref{t1}\thanksref{t2}}
 \affiliation{University of Southern California}
\address{A. G. Tartakovsky\\
University of Southern California\\
Department of Mathematics\\
3620 S. Vermont Ave, KAP-416E\\
Los Angeles, CA 90089-2532, USA\\
\printead{e2}}
\thankstext{t1}{This work was supported  by the U.S.\ Army Research Office MURI grant
W911NF-06-1-0094 and by the U.S.\ National Science Foundation grant CCF-0830419 at
the University of Southern California. The work of Moshe Pollak was also supported by a grant from the Israel Science Foundation and by the Marcy Bogen Chair of Statistics at the Hebrew University of Jerusalem.}
\thankstext{t2}{\textbf{Corresponding author.}}

\runauthor{M. Pollak and A. G. Tartakovsky}

\end{aug}

\begin{abstract}
Let $\{M_n\}_{n\ge 0}$ be a nonnegative Markov process with stationary transition probabilities. The quasistationary distributions referred to in this note are of the form
\[
\Qb_A(x) = \lim_{n\to\infty} \Pb(M_n \le x | M_0 \le A, M_1 \le A, \dots, M_{n} \le A) .
\]
Suppose that $M_0$ has distribution $\Qb_A$ and define
\[
T_A^{\Qb_A} = \min\{ n | M_n > A, n\ge 1\},
\]
the first time when $M_n$ exceeds $A$. We provide sufficient conditions for $\Eb T_A^{\Qb_A}$ to be an increasing function of $A$.
\end{abstract}



\end{frontmatter}

\section{Introduction} \label{s:Intro}

Quasistationary distributions come up naturally in the context of first-exit times of Markov processes. Of special interest --- in particular in statistical applications --- is the case of a nonnegative Markov chain, where the first time that the process exceeds a fixed level signals that some action is to be taken. The quasistationary distribution is the distribution of the state of the process if a long time has passed and yet no crossover has occurred.

Various topics pertaining to quasistationary distributions are existence, calculation, simulation, etc. For an extensive bibliography see \cite{Pollett}.

The topic addressed in this note deals with a certain aspect of the quasistationary distribution $\Qb_A$ as a function of $A$. \cite{Pollak&SiegmundJAP86} have shown, under certain conditions, that if a stationary distribution $\Qb$ exists, then $\Qb_A \to \Qb$ as $A\to\infty$.  Here we study the behavior of the expected time of the first exceedance of $A$ by a Markov process started at $\Qb_A$, as a function of $A$. Specifically, we provide conditions under which it is increasing. Our interest stems from a result in changepoint detection theory, where a certain Markov chain that calls for a declaration that  a change has taken place when a level $A$ has been exceeded has certain asymptotic optimality properties if started at the quasistationary distribution $\Qb_A$ (cf. \citealp{PollakAS85, TPP-Bernoulli2010}).

\section{Results and Examples}

Let $(\Omega, \Fc, \Pb )$ be a probability space, and let $\{M_n\}_{n=0}^\infty$ be an irreducible
Markov process defined on this space taking values in $\Mc \subseteq [0, \infty)$ and having stationary transition probabilities $\rho(t,x)= \Pb(M_{n+1} \le x | M_n =t)$.

Let $T_A = \min\set{n | M_n > A; n \ge 0}$, and assume that:
\begin{description}
\item{(C1)}  The quasistationary distribution
\[
\Qb_A(x) = \lim_{n\to\infty} \Pb(M_n \le x | T_A > n)
\]
exists for all $A > A_0 \ge 0$ (for some $A_0 < \infty$) and satisfies $\Qb_A(0) =0$.
\item{(C2)}  $\rho(s,x)$ is nonincreasing in $s$  for all fixed $x \in \Mc$.
\item{(C3)} $\rho(t s, t x)$ is nondecreasing in $t$ for all fixed $s, x \in \Mc$.
\item{(C4)} $\rho(s,  x)/\rho(s, A)$ is nonincreasing in $s$  for all fixed $x \in \Mc, x \le A $.
\item{(C5)} $\rho(t s , t x)/\rho(ts, tA)$ is nondecreasing in $t$  for all fixed $s, x, \in \Mc, x \le A$.
\end{description}

Now regard the case where $M_0$ has distribution $\Qb_A$ and define
\[
T_A^{\Qb_A} = \min\{ n | M_n > A;  n\ge 1; M_0 \sim \Qb_A \}.
\]

\begin{theorem} 
Let the  conditions {\rm (C1)}--{\rm (C5)} be satisfied. Then\\
{\rm (i)} $\Qb_{y A}(yx) \ge \Qb_A(x)$ for all $y \ge 1$ and all fixed $x \in \Mc, x \le A$;\\
{\rm (ii)} $\Eb T_A^{\Qb_A} \le  \Eb T_{yA}^{\Qb_{yA}}$ for all $y \ge 1$.
\end{theorem}

Before proving the theorem, we provide examples that show that although the conditions (C1)--(C5) are restrictive, nevertheless they are satisfied in a number of interesting cases.

Suppose $\{M_n\}_{n \ge 0}$ obeys a recursion of the form
\[
M_{n+1} = \varphi(M_n) \cdot \Lambda_{n+1}, \quad n =0, 1, \dots \,  ,
\]
where
\begin{description}
\item{(D1)}  $\{\Lambda_i\}_{i \ge 1}$ are iid positive and continuous random variables;
\item{(D2)}  the distribution function $F$ of $\Lambda_i$ satisfies
\[
\frac{F(tx)}{F(t A)} ~ \text{increases in $t$ for fixed $x\in \Mc, x \le A$} ;
\]
\item{(D3)} $\varphi(t)$ is continuous, positive and nondecreasing in $t$;
\item{(D4)} $t/\varphi(t)$ is nondecreasing in $t$;
\item{(D5)} $\varphi$ and $F$ are such that $\Pb (\lim\limits_{n\to\infty} M_n=0)=0$.
\end{description}

In this example,
\[
\rho(s,x) = F\brc{\frac{x}{\varphi(s)}}.
\]
Under these conditions, Theorem III.10.1 of \cite{Harris63} can be applied to obtain existence of a quasistationary distribution. The conditions (D1)--(D5) are easily seen to imply the conditions (C1)--(C5).

Condition (D2) is satisfied, for example, if the distribution function of $\log(\Lambda_1)$ is concave.

Many ``popular" Markov processes fit this model, some of which we now outline.\\[12pt]  
(I)  The exponentially weighted moving average (EWMA) processes:
\[
Y_{n+1} = \alpha Y_n + \xi_{n+1}, \quad n \ge 0,
\]
where $0 \le \alpha < 1$ and $\{\xi_i\}$ are iid random variables. Define $M_n = e^{Y_n}$, $\Lambda_n = e^{\xi_n}$. Here $\varphi(t) = t^\alpha$.\\[12pt]
(II)  Let $a>0$ and $\varphi(t) = t+a$, so that $M_{n+1}=(M_n +a) \Lambda_{n+1}$. When $a=1$ and $\Lambda_{n+1}$ is a likelihood ratio ($\Lambda_{n+1}=f_1(X_{n+1})/f_0(X_{n+1})$ where $X_i$ are iid), $\{M_n\}_{n\ge 0}$ is a sequence of Shiryaev-Roberts statistics for detecting a change in distribution of $X_i$, from density $f_0$ to $f_1$. The standard Shiryaev-Roberts procedure calls for setting $M_0=0$, specifying a threshold $A$ and declaring at $T_A=\min\{ n| M_n > A\}$ that a change took place. A procedure $T_A^{\Qb_A}$ that starts at a random point $M_0 \sim \Qb_A$ has asymptotic optimality properties (cf. \citealp{PollakAS85, Moustakidesetal-SS10, TPP-Bernoulli2010}). Another setting is where $r_i$ is the return on (one unit of) investment in the $i$th period and $\Lambda_i=1+r_i$, so that an investment of $m$ units at the beginning of the $i$th period will be worth $m\Lambda_i$ at its end. If one invests $a$ units at the beginning of the first period, reinvests the $a\Lambda_i$ units and adds another $a$ units at the beginning of the second period, and continues this way (i.e., always reinvesting and adding $a$ units at every period), then the process $M_{n+1} = \varphi(M_n) \Lambda_{n+1}$ with $\varphi(t) = t+a$ describes the scheme.\\[12pt]
(III) The random walk reflected from the zero barrier:
\[
Y_0=0, \quad Y_{n+1} = (Y_n+Z_{n+1})^+,  \quad n=0, 1, \dots ,
\]
where $\{Z_i\}$ are iid, $\Pb(Z_i < 0) > 0$.  Note that on the positive half plane the trajectory of the reflected random walk $\{Y_n\}_{n\ge 0}$  is identical to the trajectory of the Markov process $\{Y^{\ast}_n\}_{n\ge 0}$ given by the recursion
\[
Y_0^\ast=0, \quad Y_{n+1}^\ast  = (Y_n^*)^++Z_{n+1},  \quad n=0, 1, \dots 
\]
Therefore, if $\log A>0$ one may operate with $Y_n^\ast$ instead of $Y_n$ and all conclusions will be the same. Define $M_n = e^{Y_n^\ast}$ and $\Lambda_i=e^{Z_i}$, so that 
\[
M_{n+1} = \max(M_n, 1) \Lambda_{n+1}, \quad n \ge 0. 
\]
Here $\varphi(t) = \max(1,t)$. This process describes a broad class of single-channel queuing systems (see, e.g., \citealp{Borovkov-book76}). This setting can also be applied to the Cusum scheme  for detecting a change in distribution, when $Z_i = \log [f_1(X_i)/f_0(X_i)]$ and $X_i$, $f_0$ and $f_1$ are as in (II).


\proof
Let $\{U_n\}_{n\ge 0}$ be a Markov process with stationary transition probabilities
\[
\Pb(U_{n+1} \le x | U_n =t) = \frac{\rho(t,x)}{\rho(t, A)}, \quad x \le A,
\]
where $A>0$ is fixed and $U_0$ has an arbitrary distribution (possibly degenerate) on $[0,A]$. Let $y >1$ and define $W_n = y U_n$.

Let $\{V_n\}_{n\ge 0}$ be a Markov process with $V_0=W_0=yU_0$, having stationary transition probabilities
\[
\Pb(V_{n+1} \le x | V_n =t) = \frac{\rho(t,x)}{\rho(t, yA)}, \quad x \le yA .
\]
Clearly, the stationary distribution of $\{V_n\}$ is $\Qb_{yA}(x)$ and that of $\{W_n\}$ is $\Qb_A(x/y)$.

Since
\[
\begin{aligned}
\Pb(V_1 \le x|V_0)&  = \frac{\rho(V_0,x)}{\rho(V_0, yA)} \ge  \frac{\rho\brc{\tfrac{1}{y}V_0,\tfrac{1}{y}x}}{\rho\brc{\tfrac{1}{y}V_0, A}}\\
& = \Pb\brc{U_1 \le \tfrac{1}{y}x| U_0 = \tfrac{1}{y}V_0} = \Pb(W_1 \le x | W_0 = V_0),
\end{aligned}
\]
it follows that $V_1 \overset{\rm st} \prec W_1$ (stochastically smaller). Therefore, one can construct a sample space on which $U_0, U_1, V_0, V_1, W_0, W_1$ are all defined and such that $V_1 \ge W_1$ a.s. Write $V_1=s, W_1=t$ where $s \le t \le y A$, $s, t\in \Mc$.  Now
\[
\begin{aligned}
\Pb(V_2 \le x|V_1=s)&  = \frac{\rho(s,x)}{\rho(s, yA)} \ge  \frac{\rho(t,x)}{\rho(t, yA)}  \ge  \frac{\rho\brc{\tfrac{1}{y}t,\tfrac{1}{y}x}}{\rho\brc{\tfrac{1}{y}t, A}}\\
& = \Pb\brc{U_2 \le \tfrac{1}{y}x| U_1 = \tfrac{1}{y}t} = \Pb(W_2 \le x | W_1 = t),
\end{aligned}
\]
so that $V_2 \overset{\rm st} \prec W_2$, and one can  construct a sample space on which $U_0$, $U_1$, $U_2$, $V_0$, $V_1$, $V_2$, $W_0$, $W_1$, $W_2$ are all defined and  $V_0=W_0, V_1 \ge W_1, V_2 \le W_2$ a.s.

Continuing this inductively, one obtains a sample space on which $\{U_n\}$, $\{V_n\}$, $\{W_n\}$ are all defined and $V_n \le W_n$ a.s.\ for all $n \ge 0$. Consequently, $\lim\limits_{n\to\infty} \Pb(V_n > x) \le  \lim\limits_{n\to\infty} \Pb(W_n > x)$, i.e.,
$\Qb_{yA}(yx) \ge \Qb_A(x)$, accounting for~(i). 

To prove (ii), note that both first exit times $T_A^{\Qb_A}$ and $T_{yA}^{\Qb_{yA}}$ are geometrically distributed random variables, so that
\[
\Eb T_A^{\Qb_A}  = \frac{1}{1-\int_0^A\rho(s,A)\,d\Qb_A(s)}
\]
and
\[
\Eb T_{yA}^{\Qb_{yA}}  = \frac{1}{1-\int_0^{yA}\rho(s,yA)\,d\Qb_{yA}(s)}.
\]
Hence, it suffices to show that
\[
\int_0^{yA}\rho(s,yA)\,d\Qb_{yA}(s) \ge \int_0^{A}\rho(s,A)\,d\Qb_{A}(s) \quad \text{for $y \ge 1$}.
\]

Note that $\rho(ds, t) \le 0$. Therefore, integrating by parts yields
 \[
\begin{aligned}
& \int_0^{yA}\rho(s,yA)\,d\Qb_{yA}(s)  = \rho(s,yA)\Qb_{yA}(s) \Big |_0^{yA}-  \int_0^{yA} \Qb_{yA}(s) \rho(d s, yA)\\
 & \quad = \rho(yA, yA) - \int_0^{yA} \Qb_{yA}(s) \rho(d s, yA) ~~~ (\text{since $\Qb_{yA}(0)=0$  by (C1)}) \\
 &\quad \ge \rho(yA, yA) - \int_0^{yA} \Qb_{A}(s/y) \rho(d s, yA) ~~~(\text{by (i)})\\
&\quad  = \rho(yt, yA) \Qb_A(t) \Big | _0^A - \int_0^{A} \Qb_{A}(t) \rho(d (yt), yA)\\
&\quad = \int_0^{A}\rho(y t, yA)\,d\Qb_{A}(t)\\
& \quad\ge  \int_0^{A}\rho(t,A)\,d\Qb_{A}(t) ~~~ (\text{by condition (C3)}) ,
\end{aligned}
\]
which completes the proof.
\endproof

\appendix

\bibliographystyle{imsart-nameyear}

\end{document}